\documentclass[10pt]{article}
\usepackage[left=1.25in, right=1.25in,
    top=1in, bottom=1in]{geometry}
\usepackage[utf8]{inputenc}
\usepackage[english]{babel}
\usepackage[dvipsnames]{xcolor}
\usepackage[normalem]{ulem}
\usepackage{tikz,microtype, natbib, fancyhdr, siunitx, graphicx, hyperref}
\graphicspath{{./figures/}}

\usepackage{bm, amssymb, amsmath, amsthm}

\usepackage[sf,bf,medium]{titlesec}

\usepackage[font=footnotesize]{caption}
\usepackage[font=footnotesize]{subcaption}

\usepackage{setspace}

\newcommand{\R}{\mathbb{R}}    
\newcommand{\bO}{\mathcal{O}}  
\newcommand{\Nd}{\mathcal{N}}  
\newcommand{\sv}{\, | \,}      
\newcommand{\bx}{\bm{x}}       
\newcommand{\by}{\bm{y}}       
\newcommand{\bv}{\bm{v}}       
\newcommand{\bu}{\bm{u}}       
\newcommand{\bS}{\bm{\Sigma}}  
\newcommand{\abs}[1]{\left|#1\right|}      
\newcommand{\norm}[2][]{\left\Vert#2\right\Vert_{#1}} 
\newcommand{\set}[1]{\left\{ #1 \right\}}             
\newcommand{\mat}[1]{\begin{bmatrix}#1\end{bmatrix}}  

\renewcommand{\Re}{\operatorname{Re}}

\newcommand{\lp}{\left(}
\newcommand{\rp}{\right)}

\newtheorem*{theorem*}{Theorem}

\newtheorem*{lemma*}{Lemma}

\newtheorem*{corollary*}{Corollary}

\newtheorem*{definition*}{Definition}

\newtheorem*{proposition*}{Proposition}

\numberwithin{equation}{section}

\newcommand{\tps}{\varphi_{\text{tps}}}
\newcommand{\Qperp}{\bm{Q}_{\perp}}

\title{\sffamily \bfseries Linear-cost Polyharmonic Spline Interpolation of
  Arbitrary Degree}
\author{Christopher J. Geoga\footnote{Department
        of Statistics, University of Wisconsin-Madison, Madison, WI, 53706.
        \texttt{geoga@wisc.edu}.}%
        \and%
        Michael O'Neil\footnote{Department of Mathematics, Courant Institute
        School of Mathematics, Computing, and Data Science, New York
        University, New York, NY 10012. \texttt{oneil@nyu.edu}.}}
\date{}
  
\begin{document}

\maketitle

\begin{abstract}

We introduce a simple and performant approach for rapidly and accurately
performing polyharmonic spline (PHS) interpolation using a combination of the
fast multipole method (FMM) from computational electrostatics and the Vecchia
approximation from the Gaussian process and sparse approximate inverse
literatures. Using basic properties about Hadamard products and low-rank
matrices, we demonstrate that an FMM with two kernels, the logarithmic and
distance kernels, results in fast PHS interpolation for \emph{all} orders.
Furthermore, we demonstrate the exceptional performance of sparse inverse
approximation methods with the Mat\'ern covariance model for preconditioning.
Combined with careful management of disallowed subspaces, we describe a
procedure for obtaining prediction weights using preconditioned conjugate
gradient that converges in less than $15$ iterations, even for problem sizes
with over one million points. As a result, thin-plate spline interpolation---a
particularly popular method that does not require parameter tuning---that
matches the fully dense $\bO(n^3)$ computation in accuracy can be done at the
cost of approximately $50-60$ FMMs.  A high-performance software library for
odd-order PHS interpolation in two dimensions is made available as a companion
to this work.

\end{abstract}


\setstretch{1.1}

\section{Introduction} \label{sec:intro}

Polyharmonic spline (PHS) interpolation is an extremely successful and popular
variety of kernel interpolation methods \cite{wahba1990}. Given data $\set{(\bx_j,
f(\bx_j))}_{j=1}^n$ PHS interpolation involves decomposing the estimation of $f$
at un-measured locations into a polynomial trend term and a dense kernel
interpolation-based remainder, so that
\begin{equation} \label{eq:phs_interp}
  \hat{f}(\bx) =  \sum_{s=1}^S \beta_s P_s(\bx) + \sum_{j=1}^n w_j
  \varphi(\norm{\bx_j - \bx}),
\end{equation}
where~$\hat{f}$ is the PHS interpolant, $\set{P_s}_{s=1}^S$ are tensor products
of low-degree monic polynomials in each variable, $\varphi$ is the polyharmonic
spline radial basis function defined as
\begin{equation} \label{eq:phs_fun}
  \varphi(r) = 
  \begin{cases}
    r^k        & \text{$k \in \mathbb{N}$ odd} \\
    r^k \log r & \text{$k \in \mathbb{N}$ even}
  \end{cases},
\end{equation}
and $\bm{w} = [w_j]_{j=1}^n$ are interpolating weights. The weights~$w_j$ and the
coefficients~$\beta_s$ are given as the solution to the augmented linear
system
\begin{equation} \label{eq:phs_linsys}
  \mat{
    \bm{M} & \bm{B} \\
    \bm{B} & \bm{0}
  }
  \mat{\bm{w} \\ \bm{\beta}}
  =
  \mat{\bm{f} \\ \bm{0}},
\end{equation}
where $\bm{M}_{j,k} = \varphi(\norm[2]{\bx_j - \bx_k})$ is the kernel matrix
generated by the PHS function $\varphi$, $\bm{B}_{j,k} = P_j(\bx_k)$ is the
matrix of polynomial functions up to (joint) order $\geq k$, and $\bm{f} =
[f(\bx_j)]_{j=1}^n$ is the collection of given measurements. This specific
formulation for PHS interpolation is necessary because the kernel function
$\varphi$ is only \emph{conditionally} positive-definite: if $\bm{B} = \bm{Q}
\bm{R}$ is the thin QR decomposition of $\bm{B}$ chosen to contain exactly
polynomials up to degree $k$, then the eigenvalues of $\bm{Q}^T \bm{M} \bm{Q}$
are negative and $\bm{Q}_{\perp}^T \bm{M} \bm{Q}_{\perp} \succ \bm{0}$ for
$\bm{Q}_{\perp}$ such that $[\bm{Q} \; \bm{Q}_{\perp}]$ is an $n \times n$
orthogonal matrix \cite{wahba1990,wendland2005}, ensuring that weights for an
interpolating solution $\hat{f}$ exist. Furthermore, unlike many other kernel
interpolation methods, PHS interpolation is scale-invariant, meaning that
lengthscale parameters do not need to be tuned to produce performant
interpolation schemes.

The challenge with many kernel interpolation methods (including PHS approaches)
is that $\bm{M}$ is a dense matrix that is prohibitively expensive to factorize,
and potentially even to form, for large data sizes $n$. Direct methods to
solve~\eqref{eq:phs_linsys} will scale with $\bO(n^3)$ complexity, as they
require a factorization of the augmented $(n + S) \times (n + S)$ matrix (with
the Bunch-Kaufman or ``LBLt" factorization \cite{bunch1977} being the standard
choice). A popular approach to addressing this is to use iterative methods
\cite{saad2003}, which if one uses exact matrix-vector products can reduce this
cost to $\bO(n^2)$, presuming that a sufficiently accurate preconditioner can be
designed to control the number of iterations required for convergence. Several
existing works in the literature offer approaches of this category. In
\cite{beatson1999}, for example, variants of the fast multipole method
(FMM)~\cite{greengard1987fast} are given for specific degrees, and a
preconditioned GMRES approach is given that ostensibly brings the cost of
solving~\eqref{eq:phs_linsys} down to~$\bO(n)$ if one uses a careful and
non-standard stopping criterion. As we will demonstrate below, however, it does not
control iteration counts under standard convergence metrics.  Similarly, in
\cite{powell1994}, a clever management of the column space of $\bm{B}$ allows
for preconditioned conjugate gradient to be used, although the reported
iteration counts again clearly grow with $n$ using the specific preconditioning
strategy proposed there.

The main contribution of this work is to simplify and improve the fast
matrix-vector products and preconditioning approaches. First, we describe a
simple modification to the matrix-vector computation $\bv \mapsto \bm{M} \bv$
that exploits a standard identity for Hadamard products with low-rank matrices.
As a result of this modification, the application $\bv \mapsto \bm{M} \bv$ is
computable for \emph{any} spline order $k$ using a single FMM routine for a
``core'' kernel of $\varphi(r) = r$ or~$\varphi(r) = \log r$.  Second, we
describe a simpler preconditioner than those proposed in prior work that only
requires~$\bO(n)$ floating point operations to apply, is easily assembled in
parallel with no required communication between workers, and is substantially
more effective than its predecessors. Put together, these modifications provide
truly end-to-end linear-cost PHS interpolation that are easy and convenient to
implement using existing well-tested and popular software. Using the
library~\texttt{fmm2d}~\cite{askham2021fmm2d} for core matrix-vector products,
for example, we provide demonstrations of $n=10^6$ points being interpolated at
$m=10^4$ new locations in under $15s$ on a modest laptop.

\section{Fast Computation of PHS Weights} \label{sec:method}
The algorithmic simplifications and enhancements of this work depend on two
critical observations about the interpolation problem. The first of these
observations is that the kernel component of PHS interpolation is \emph{scale
  invariant}. Throughout the exposition of this and subsequent sections, for
simplicity we will use the \emph{thin-plate spline} function
\begin{equation*} 
  \tps(r) = r^2 \log r.
\end{equation*}
Indeed, this function~$\tps$ is the polyharmonic spline $\varphi$ from the
introduction with $k=2$. Where appropriate (and when not obvious), we will
comment on how methods and observations can be extended to higher polynomial
orders $k>2$.

For demonstration, consider the matrix $\bm{M}_{j,k} = \tps(\norm[2]{\bx_j -
\bx_k})$. Clearly, if each location is shifted as $\bx_j \mapsto \bx_j +
\bm{\xi}$, the matrix $\bm{M}$ is unchanged.  But more importantly, if each
location is \emph{scaled} as $\bx_j \mapsto \rho \bx_j$, we note that
\begin{equation*} 
  \tps(\norm[2]{\rho(\bx_j - \bx_k)})
  =
  \rho^2 r^2 \log r + \rho^2 \log(\rho) r^2.
\end{equation*}
Letting $\bm{M}(\rho)_{j,k} = \tps(\norm[2]{\rho(\bx_j - \bx_k)})$ denote the
matrix of scaled locations, we note from above that
\begin{equation} \label{eq:M_rho}
  \bm{M}(\rho) 
  = 
  \rho^2 \bm{M} + \rho^2 \log \rho [\norm[2]{\bx_j - \bx_k}^2]_{j,k=1}^n,
\end{equation}
where
\begin{equation}
  [\norm[2]{\bx_j - \bx_k}^2]_{j,k=1}^n
\end{equation}
denotes the~$n \times n$ matrix whose~$j,k$ entry is given by~$\norm[2]{\bx_j - \bx_k}^2$.
However, observing that 
\begin{equation} \label{eq:norm_sep}
  \norm[2]{\bx_j - \bx_k}^2 = \norm[2]{\bx_j}^2 + \norm[2]{\bx_k}^2 - 2\bx_j^T
  \bx_k,
\end{equation}
we see that the second term in~\eqref{eq:M_rho} is exactly comprised of quadratic
polynomials of the first, second, and mixed components of $\set{\bx}_{j=1}^n$.
These terms, and thus polynomial updates to the interpolation points $\bm{f}$,
do not impact $\bm{w}$ due to the constraint $\bm{B} \bm{w} = \bm{0}$ introduced
in the second row of (\ref{eq:phs_linsys}). Let
\begin{equation}
  \bm{B} = [P_k(\bx_j)]_{j=1,k=1}^{n,S} = \bm{Q} \bm{R}
\end{equation}
be the thin truncated QR of the polynomial span to be removed, and let $\Qperp$
be a matrix with orthogonal columns such that $[\bm{Q} \; \Qperp]$ is an $n
\times n$ orthogonal matrix. By construction, then, $\Qperp^T \bm{M}(\rho)
\Qperp = \rho^2 \Qperp^T \bm{M} \Qperp$. And since~\eqref{eq:phs_linsys} can
practically be solved by sequentially solving the systems
\begin{align} \label{eq:phs_linsys_solve}
  \Qperp^T \bm{M} \Qperp \bm{u}
  &= \Qperp^T \bm{f}
  \\
  \notag
  \bm{w} &= \Qperp \bm{u}
  \\
  \notag
  \bm{R} \bm{\beta} &= \bm{Q}^T(\bm{f} - \bm{M} \bm{w}),
\end{align}
we see that rescaling the locations does not impact the dense $n \times n$
linear system at all. As such, we are free to rescale points to a domain that is
convenient for our method without impacting the bottleneck of obtaining
interpolation weights $\bm{w}$. This formulation using $\Qperp$ was first used
in \cite{powell1994}, who observed that $\Qperp$ can be applied to a matrix at a
cost of $\bO(n S)$ by exploiting the fact that
\begin{equation*} 
  \bm{H} \bu
  =
  \begin{bmatrix}
    \bm{Q}^T \bu \\
    \Qperp^T \bu
  \end{bmatrix},
  \quad
  \quad
  \bm{H} = \prod_{j=1}^S (\bm{I} - 2 \bv_j \bv_j^T),
\end{equation*}
where $\bm{H}$ is precisely the matrix that is formed from the thin Householder
QR factorization of $\bm{B}$, since there are only $S$ reflectors and each can
be applied in $\bO(n)$ time.  This simple trick transforms the dense and
indefinite $n \times n$ linear sub-system of~\eqref{eq:phs_linsys} into an $(n -
S) \times (n - S)$ system that is strictly positive-definite. Thus, PCG can be
used to solve for $\bm{u}$ (and thus obtain weights $\bm{w}$). A key difference
between this proposed method and that of \cite{powell1994}, however, is how
structure in this formulation will be used to build an effective preconditioner
for $\bm{M}$, as it need only be accurate for the positive definite part of
$\bm{M}$ that is not annihilated by $\Qperp$. 

The second powerful and related observation pertains to the application of the
matrix~$\bm{M}$. As mentioned above, we note that the matrix generated by powers
$r^k$ for even $k$ is rank deficient because norms raised to even powers are
separable in their arguments.  Specifically, in two dimensions, we see that
\begin{equation*} 
  [\norm[2]{\bx_i - \bx_j}^k]_{i,j=1}^n
  =
  \bm{U} \bm{V}^T,
\end{equation*}
where $\bm{V} \in \R^{n \times 2^{2(k-1)}}$. In the case of $k=2$ in two
dimensions, for example, we see that these matrices are given by 
\begin{align*} 
   \bm{U}_{j,\cdot}
   &=
   \begin{bmatrix}
     \norm[2]{\bx_j}^2 & 1 & -2 x_j & -2 y_j
   \end{bmatrix},
   \\
  \bm{V}_{j,\cdot}
  &=
  \begin{bmatrix}
    1 & \; \; \; \norm[2]{\bx_j}^2 & \; \; x_j & \; \; \; \; y_j
  \end{bmatrix},
\end{align*}
as expected, given the decomposition~\eqref{eq:norm_sep}. As a result, we
observe that $\bm{M} = \bm{M}_{\text{log}} \circ (\bm{U} \bm{V}^T)$, where
$[\bm{M}_{\text{log}}]_{j,k} = \log r$ for $j \neq k$ and $0$ on the diagonal is
what we refer to as a ``core" kernel matrix. Using standard properties for
Hadamard products with low rank matrices, we see that
\begin{equation} \label{eq:matvec_hadamard}
  \bm{M} \bm{t} = \sum_{l=1}^{2^{2(k-1)}} \bm{D}_{\bu_l}
  \bm{M}_{\text{log}} \bm{D}_{\bv_l} \bm{t},
\end{equation}
where $\bm{D}_{\bu_l} = \text{Diag}(\bm{U}_{\cdot,l})$ and
$\bm{D}_{\bv_l} = \text{Diag}(\bm{V}_{\cdot,l})$. Due to this structure, we see
that for \emph{any} order~$k$, the problem of rapidly computing the
matrix-vector product $\bm{u} \mapsto \bm{M} \bm{u}$ effectively reduces to a
fast matrix-vector product with $\bm{M}_{\text{log}}$. Fortunately, the
kernel~$\log r$ is the Green's function for the Poisson equation in two
dimensions, and extremely high-performance computational algorithms for the
associated~$n$-body sums have been developed over the past several decades. We
now briefly describe the idea behind these algorithms that sit at the core of
this simplified approach.

\subsection{The Fast Multipole Method}
\label{subsec:fmm}

As discussed above, variants of the classical version of the fast multipole
method (FMM) have been constructed to accelerate the PHS interpolation problem
in two dimensions in \cite{beatson1999}. The previous discussion makes the
\textit{new} observation that all that is necessary, in fact, for computing the
matrix-vector product~$\bm{M} \bv$ for \emph{any} order PHS are matrix-vector
products with ``core'' kernels $\Vert \bm{r}\Vert$ and~$\log\Vert\bm{r}\Vert$.
(For convenience in the following section, will we generally use the notation~$r
= \Vert \bm{r} \Vert$, where~$\Vert \cdot \Vert$ denotes the usual Euclidean
distance in two dimensions.) In this section, we will give a high-level
description of the linear-cost FMM algorithm for the $\log$ kernel.  The FMM
literature is deep and vast, so we will not give an exhaustive list of
references, but merely point the reader to a couple specific works that are of
particular relevance for our application~\cite{beatson1997short,
greengard1987fast}.

Application of the matrix~$\bm{M}_{\text{log}}$
in the previous section can be boiled down to the evaluation of an~$N$-body sum
\begin{equation}
  \label{eq:nbody}
  \psi(\bm{x}_i) = \sum_{j=1}^{N} w_j \, \log\Vert \bm{x}_i - \bm{y}_j\Vert,
\end{equation}
for~$i = 1,\ldots,M$, where we have changed notation slightly to clearly mark
what we will refer to as \textit{targets}~$\bm{x}_i$ and
\textit{sources}~$\bm{y}_j$. In many cases, the set of targets is the same as
the set of sources (as in the previous discussion) and therefore we must
eliminate the~$j=i$ term in the sum.  At its core, the FMM makes use of the fact
that the kernel function~$\log r = \log\Vert\bm{x} - \bm{y}\Vert$ can be
approximated, to an arbitrarily high precision, via an expansion which is
separable in the variables~$\bm{x}$ and~$\bm{y}$ (assuming some distance of
separation between the points). This means we can write
\begin{equation}
  \log\Vert \bm{x} - \bm{y} \Vert \approx \sum_{k=0}^p u_k(\bx) \, v_k(\by), 
\end{equation}
for some known functions~$u_k$,~$v_k$.
Systematically grouping the sum in~\eqref{eq:nbody} depending on the locations
of the sources and targets, and using a separable expansion of the~$\log$ kernel,
allows for the design of an algorithm which can compute the above sum
in~$\mathcal O(M + N)$ time, accurate to any user-specified precision.

To this end, let us briefly examine the exact separable expansion for the~$\log$
kernel. Temporarily, let~$w$ and~$z$ be complex numbers such that~$|w|>|z|$.
We begin with the standard Taylor series expansion
for~$\log$:
\begin{equation}
  \log(w - z) = \log w - \sum_{k=1}^\infty \frac{1}{k} \left( \frac{z}{w}
  \right)^k.
\end{equation}
Also note that if we associate the point~$\bm{x} \in \mathbb R^2$ with~$w \in
\mathbb C$ and~$\bm{y} \in \mathbb R^2$ with~$z \in \mathbb C$, then we have
that
\begin{equation}
  \log\Vert \bm{x} - \bm{y} \Vert = \Re\left(\log(w-z) \right).
\end{equation}
Given the above relationship between the complex-valued logarithm and the
real-valued one, it suffices to describe the relevant expansions and algorithm
in the context of the complex-valued potential and then take the real-part (this
simplifies both the algebra and numerical implementations as well -- the
software library \texttt{fmm2d} is based on such an approach~\cite{askham2021fmm2d}).

The above expression for~$\log(w-z)$ can be truncated after~$p$ terms and
bounded by an error~$\epsilon$,
\begin{equation}
  \log(w - z) = \log w - \sum_{k=1}^p \frac{1}{k} \left( \frac{z}{w}
  \right)^k + \epsilon, 
\end{equation}
where
\begin{equation}
  |\epsilon| < \frac{1}{p+1} \sum_{k=p+1}^\infty \left| \frac{z}{w}
  \right|^{p+1} = \left( \frac{1}{p+1} \right) \left( \frac{1}{|z/w| - 1}
  \right)
  \left| \frac{z}{w} \right|^p.
\end{equation}
In particular, if~$|w| > 2|z|$, then~$|\epsilon| < 1/2^p$. This gives a rigorous
bound for the above expansion of the~$\log$ kernel. Indeed, if the
source~$\bm{y}$ and target~$\bm{x}$ are separated by more than a factor of 2,
relative to the origin,~$p$ can be chosen using this error bound to ensure any
specified accuracy~$\epsilon$. In this regime, the expansion is accurate to
about~$10^{-6}$ with~$p=20$ and $10^{-16}$ with~$p=50$.

Next, consider a collection of~$N$ sources~$z_j$ all contained within a ball of
radius~$R$, i.e.~$|z_j|<R$, each with strength~$q_j$. Defining a potential~$\Phi$
at a collection of~$M$ targets~$w_i$, with~$|w_i|=r>R$, by the sum
\begin{equation}
  \Phi(w_i) = \sum_{j=1}^N q_j \, \log(w_i-z_j),
\end{equation}
we can invoke the above approximation to yield a ``fast'' algorithm for
evaluating the~$N$-body sum (i.e. the matrix-vector problem with the~$\log$
kernel).
We have
\begin{equation}
  \label{eq:outgoing}
  \begin{aligned}
    \Phi(w_i) &= \sum_{j=1}^N q_j \, \log(w_i-z_j) \\
              &= \sum_{j=1}^N q_j \left( \log w_i - \sum_{k=1}^\infty
                \frac{1}{k} \left( \frac{z_j}{w_i}
                \right)^k \right)  \\
              &= Q \log w_i - \sum_{k=1}^\infty \lp
                \sum _{j=1}^N q_j \frac{z_j^k}{k}
                \rp
                \frac{1}{w_i^k} \\
              &= Q \log w_i - \sum_{k=1}^\infty a_k 
                  \frac{1}{w_i^k}.
  \end{aligned}
\end{equation}
Of course the above sum could be truncated after~$p$ terms, with~$p$ determined
similarly as before, making the finite expansion accurate to an arbitrary
precision. The above expansion, given in terms of the powers~$1/w_i^k$ will be
referred to as an \textit{outgoing} expansion. By simple counting, the above
potential can be evaluated at all of the targets~$w_i$ using~$\mathcal
O(p(M+N))$ operations.

Conversely, if the role of targets and sources were switched, we could write
down a similar expansion in terms of powers~$z_j$:
\begin{equation}
  \label{eq:incoming}
  \begin{aligned}
    \Phi(z_j) &= \sum_{i=1}^M \tilde q_i \, \log(w_i-z_j) \\
              &= \sum_{i=1}^M \tilde q_j \left( \log w_i - \sum_{k=1}^\infty
                \frac{1}{k} \left( \frac{z_j}{w_i}
                \right)^k \right)  \\
              &= \sum_{i=1}^M \tilde q_i \log w_i - \sum_{k=1}^\infty \lp
                \sum _{i=1}^M \tilde q_i \frac{1}{k w_i^k}
                \rp
                z_j^k \\
              &=  \sum_{k=1}^\infty b_k z_j^k.
  \end{aligned}
\end{equation}
Expansions of this form will be referred to as \emph{incoming expansions}, and
can also be rapidly evaluated at all the points~$z_j$ at a cost
of~$\mathcal{O}(p(M+N))$.  The main difficulty lies when the sources and targets
are interlaced (i.e. not geometrically separated) and therefore a single
outgoing (or incoming) expansion cannot be used. In this regime, a hierarchical
data structure must be used to systematically sort sources and targets into
regions which are well-separated and whose potential interactions then admit
outgoing and incoming expansions which are accurate to high precision using
$p$-term expansions with a modestly sized \emph{fixed} value of~$p$.

\begin{figure}[t!]
  \centering
  \begin{subfigure}[b]{0.45\linewidth}
    \centering
    \includegraphics[width=.95\linewidth]{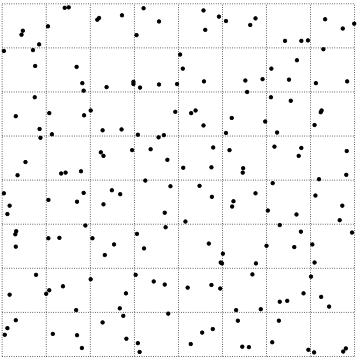}
    \caption{Source-target configuration.}
  \label{fig:fmm1}
  \end{subfigure}
  \hfill
  \begin{subfigure}[b]{0.45\linewidth}
    \centering
    \includegraphics[width=.95\linewidth]{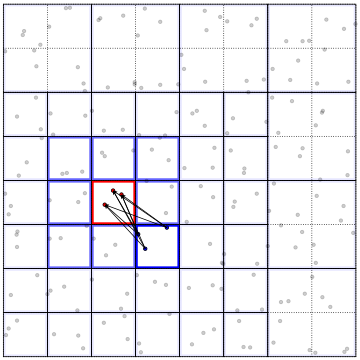}
    \caption{Direct evaluations.}
  \label{fig:fmm2}
  \end{subfigure}\\
  \vspace{\baselineskip}
  \begin{subfigure}[b]{0.45\linewidth}
    \centering
    \includegraphics[width=.95\linewidth]{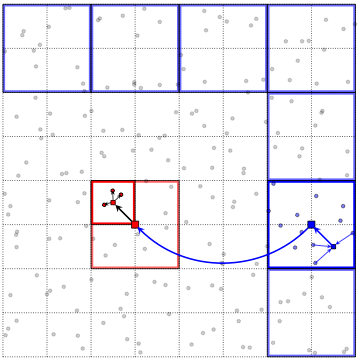}
    \caption{Top-level field translations.}
  \label{fig:fmm3}
  \end{subfigure}
  \hfill
  \begin{subfigure}[b]{0.45\linewidth}
    \centering
    \includegraphics[width=.95\linewidth]{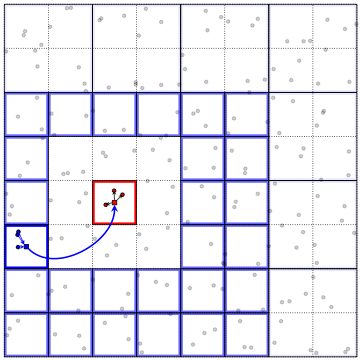}
    \caption{Leaf-level field translations.}
  \label{fig:fmm4}
  \end{subfigure}
  \caption{Geometric depiction of a fast multipole method. Reproduced with
    permission from the authors from the documentation
    of \texttt{FMM3D}~\cite{fmm3d}.}
  \label{fig:fmm}
\end{figure}

To precisely describe the full algorithm would be beyond the scope of this work,
as both fundamentals and implementation specifics of fast multipole methods is
well-established at this point in time. However, we will provide a graphical
overview of the algorithm to give the reader an idea of how it works. Consider the
point configurations and groupings in Figure~\ref{fig:fmm}.
The general fast multipole method proceeds via the following steps:
\begin{enumerate}
\item \emph{Tree creation.} All sources and targets are sorted into a quadtree
  data structure, with leaf boxes containing at most some fixed number of
  points, see Figure~\ref{fig:fmm1}.
\item \emph{Direct evaluation.}  For leaf boxes which are adjacent to each
  other, i.e. share a corner (including itself), the potential due to the
  sources within is directly evaluated at each of the targets in the self
  box. See Figure~\ref{fig:fmm2}.
\item \emph{Upward pass.} For each leaf box, an outgoing expansion of the
  form~\eqref{eq:outgoing} is created which represents the potential due to all
  sources in this box, valid at all targets in \emph{well-separated} boxes,
  i.e. those separated from the source box by at least another box. Similar
  outgoing expansions for subsequent \emph{parent} boxes can be constructed from
  the outgoing expansions for \emph{children} boxes via a linear operator known
  as the \emph{outgoing-to-outgoing translation operator}. See blue arrows in
  blue box on the right in Figure~\ref{fig:fmm3}. This proceeds to the coarse
  level shown in Figure~\ref{fig:fmm3}.
\item \emph{Outgoing-to-incoming translations} At this point, all boxes on all
  levels in the quadtree hierarchy have associated with them an outgoing
  expansion of the form~\eqref{eq:outgoing} representing the potential due to
  all sources contained within, and valid at all targets contained in
  well-separated boxes. Each of these outgoing expansions is now transformed to an
  incoming expansion of the form~\eqref{eq:incoming} via a linear
  \emph{outgoing-to-incoming translation operator}. Each box \emph{receives} a
  translation from every box which is a child of its parent's neighbors, and
  which is not adjacent (this list is known as the interaction list for a
  box). See Figures~\ref{fig:fmm3} and~\ref{fig:fmm4} for depictions of
  interaction lists for a red box; the blue arrow connecting a red box and a
  blue box represents application of the outgoing-to-incoming translation
  operator.
\item \emph{Downward pass.} All that remains now is to push each box's incoming
  expansion---which represents the potential due to all well-separated
  sources---down to the leaf level. This can be accomplished by starting at the
  coarsest level of the quadtree hierarchy, and for each box, apply an
  \emph{incoming-to-incoming translation operator} which transfers the information
  encoded in the incoming expansion of a parent box to an incoming expansion of
  its children. These translations are depicted using the red arrows in
  Figure~\ref{fig:fmm3}.
\item \emph{Local evaluations.} After the previous downward pass, each box on
  the leaf level now contains an incoming expansion which represents the
  potential due to \emph{all} well-separated boxes, i.e. exactly those boxes
  whose contribution was not accounted for in Step 2, the direct evaluation
  step. Each of these local incoming expansions is then evaluated at each of the
  targets inside a leaf box, and the potential is added to the previous
  potential computed from the neighboring boxes. See the red arrows inside the
  leaf box in Figure~\ref{fig:fmm4}.
\end{enumerate}

It should be clear from the above steps that the work performed for each box is
purely local, and comes from either a direct potential evaluation or application
of one of the translation operators. The translation operators map~$p$
coefficients in one expansion (outgoing or incoming) to~$p$ coefficients in
another expansion (outgoing or incoming) and can be shown to be dense
operators. Therefore their naive application requires~$\mathcal O(p^2)$
work. Assuming there are~$N$ sources and targets, it can be shown that there
are~$\mathcal O(N)$ boxes, and the overall cost for the previous algorithm is
therefore~$\mathcal O(p^2 N)$.

Of course, there are many variants of the above basic algorithm which use
adaptive data structures and parallel computing paradigms, but the basic idea is
the same. We did not describe the exact form of the various translation
operators mentioned in the above cursory description of the algorithm, but they
can be obtained via relatively simple algebraic manipulations of the outgoing
and incoming expansions, and are given in~\cite{greengard1987fast}. Various
accelerations of the application of these translation operators can be derived
based on alternative representations of the outgoing and incoming potentials
known as \emph{plane-wave} or \emph{exponential}
representations~\cite{greengard1998}. These accelerated translation operators
are implemented in the \texttt{fmm2d} library on which we base our algorithm.
They reduce the prefactor in the asymptotic cost, but not the asymptotic cost as
a function of~$N$.

\subsection{Preconditioning with Vecchia's Approximation and Mat\'ern kernels}
\label{subsec:vecchia}

In this section, we propose a high-performance and extremely fast $\bO(n)$-cost
sparse preconditioner for the thin-plate spline matrix $\bm{M}$ defined above.
To begin, we first note that the thin-plate spline function $\tps$ is one of the
two linearly independent fundamental solutions of the \emph{biharmonic equation}
\begin{equation} \label{eq:biharmonic}
  \Delta^2 \tilde{\varphi} = 0
\end{equation}
on $\R^2$ \cite{wahba1990}. This motivates a connection with the
Whittle-Mat\'ern stochastic partial differential equation, which describes a
Gaussian process $Z$ on~$\R^2$ given by the solution to
\begin{equation*} 
  (\alpha^2 - \Delta)^{\alpha/2} Z = \mathcal{W}
\end{equation*}
where $\mathcal{W}$ is spatial white noise in~$\R^2$~\cite{whittle1963}. It
is well known that the covariance between $Z(\bx)$ and $Z(\bx')$ is,
up to proportionality, given by
\begin{equation*} 
  \text{Cov}(Z(\bx), Z(\bx'))
  =
  K(\bx - \bx')
  \propto 
  (\alpha \norm[2]{\bx - \bx'})^{\nu} \,  \mathcal{K}_{\nu}(\alpha \norm[2]{\bx - \bx'}),
\end{equation*}
where $\mathcal{K}_{\nu}$ is the modified second-kind Bessel function of order
$\nu$ \cite{NIST}, and $\alpha = \nu + 1$ in two dimensions. In the spatial
statistics literature, $\nu$ is called the \emph{smoothness} parameter, as it
dictates the degree of (potentially fractional) mean-square differentiability of
the process~$Z$~\cite{matern1960,whittle1963,stein1999}.  Setting~$\nu=1$ and
observing that (\ref{eq:biharmonic}) describes a transport map, we see that the
corresponding precision operator is $(\alpha^2 - \Delta)^2 = \Delta^2 - 2
\alpha^2 \Delta + \alpha^4$, from which we see the~$\Delta^2$ term from the
biharmonic equation appear.  By the expansion of the precision operator, one
would expect high-frequency behavior of the Mat\'ern and $\tps$ functions to
agree well. Taking a Taylor expansion of the Mat\'ern kernel for small arguments
$r = \norm[2]{\bx - \bx'}$ confirms this, as
\begin{equation} \label{eq:matern_pre}
  \alpha r \mathcal{K}_{1}(\alpha r)
  =
  \underbrace{1 + c_M (\alpha r)^2}_{\text{polynomial terms}}
  +
  \underbrace{\frac{(\alpha r)^2}{2} \log \frac{\alpha r}{2}}_{\propto \phi(r)}
  + 
  \underbrace{\bO((\alpha r)^4 \log (\alpha r))}_{\text{small when $\alpha r < 1$}}.
\end{equation}
As the underbrace labels indicate, for nearby points we see that the leading
term of $K(\bx - \bx')$ matches $\tps$ up to scaling and the polynomial
contributions that are annihilated by the $\Qperp$ projection in
(\ref{eq:phs_linsys_solve}).  With this in mind, we propose preconditioning
$\Qperp^T \bm{M} \Qperp$ with $\Qperp^T \bS^{-1} \Qperp$, with $\bS = [K(\bx_j -
\bx_k)]_{j,k=1}^n$ with $\nu=1$ and~$\alpha$ selected based an on error balancing
procedure discussed later in this section. And while we note that of course
$(\bm{Q}_{\perp}^T \bm{M} \bm{Q}_{\perp})^{-1} \neq \bm{Q}_{\perp}^T \bm{M}^{-1}
\bm{Q}_{\perp}$, since the matrix~$\bm{Q}_{\perp}$ is in~$\R^{n \times (n-S)}$,
for small $S$, ignoring the Schur complement term in the inverse will only add a
small number of iterations to an iterative solver routine. In our
experimentation, this provides faster start-to-finish runtimes than accounting
for the low-rank update.  As an additional complication, $\bS$ is a dense $n
\times n$ matrix itself, and so direct factorization and linear system solving
is prohibitively expensive. We will now discuss the remedy for this
computational infeasibility using tools from the Gaussian process literature.
Following that, we will discuss methods to sharpen the approximation based
on the scale-invariance of the PHS functions $\varphi(r)$ in the column space of
$\Qperp$.

\emph{Vecchia}'s approximation
\cite{vecchia1988,stein2004,finley2019,katzfuss2021, schaefer2021sparse} for a
Gaussian process log-likelihood (also known as a \emph{nearest neighbor Gaussian
process} approximation or a \emph{factorized sparse approximate inverse} (FSAI)
in the applied math community \cite{kaporin1994,yeremin2000}) is based on the
very simple idea of approximating a sequence of conditional distributions. An
elementary result of probability is that joint probabilities can be approximated
with telescoped conditionals (via repeated applications of $p(x, x') = p(x' \sv
x) p(x)$).  Presume, for example, that we aim to approximate the log-density for
$\by \sim \Nd(\bm{0}, \bS)$, where $\bS$ is given by
\begin{equation}
  \bS_{j,k} = \text{Cov}(Y(\bx_j), Y(\bx_k))
  = K(\bx_j - \bx_k).
\end{equation}
The associated joint log-density can be expanded as
\begin{equation*} 
  \log p(\by) = \log p(y_1) + \sum_{j=2}^n \log p(y_j \sv y_1, ..., y_{j-1}).
\end{equation*}
We also note that these conditionals have a simple concrete form, as
\begin{align} \label{eq:conditionals}
  y_j \sv y_1, ..., y_{j-1}
  &\sim 
  \Nd\left(\bm{\lambda}_j^T \by_{j-1}, 
  \; \;
  K(\bx_j, \bx_j) - 
  K\left(\bx_j, \underline{\bm{X}}_{j-1} \right)^T \bm{\lambda}_j\right),
  \\
  \notag
  \bm{\lambda}_j &= 
  K\left(\underline{\bm{X}}_{j-1},
  \underline{\bm{X}}_{j-1}\right)^{-1}
  K(\bx_j, \underline{\bm{X}}_{j-1}),
\end{align}
where $\underline{\bm{X}}_{j-1} = \set{\bx_k}_{k=1}^{j-1}$, $K(\bx_j,
\underline{\bm{X}}_{j-1}) = [K(\bx_j - \bx_k)]_{k=1}^{j-1}$, and
$K\left(\underline{\bm{X}}_{j-1}, \underline{\bm{X}}_{j-1}\right)$ is the
marginal covariance matrix for the process at locations
$\underline{\bm{X}}_{j-1}$. These $\bm{\lambda}_j$ weights are the \emph{best
linear prediction} weights, or \emph{Kriging} weights. For more background on
the theory of linear prediction for Gaussian processes, we refer the reader to
\cite{stein1999}.  The idea of Vecchia's approximation in this context is very
natural: if the process~$Z$ has Markovian-like structure (referred to as the
\emph{screening} property in the GP literature \cite{stein2011}), then the
conditional distribution $y_j \sv y_1, ..., y_{j-1}$ can be well-approximated
with $y_j \sv y_{\sigma(j)_1}, ..., y_{\sigma(j)_k}$, where $\sigma(j) \subset
[j-1]$ is a set of size $|\sigma(j)|=k$ called the \emph{conditioning set}.
Oftentimes, one chooses conditioning sets to contain indices of prior enumerated
points that correspond to measurements at locations that are close in physical
space to~$\bx_j$. From this idea, one can immediately obtain a good
approximation for the joint log-density of $\by$ with
\begin{equation}
  \label{eq:vecchia_sumform}
  \begin{aligned}
  \log p(\by) 
  &= 
    \log p(y_1) + \sum_{j=2}^n \log p(y_j \sv y_1, ..., y_{j-1})\\
  &\approx
    \log p(y_1) + \sum_{j=2}^n \log p(y_j \sv y_{\sigma(j)}).
  \end{aligned}
\end{equation}
If $|\sigma(j)| = \bO(1)$ for all $j$, then this immediately gives a trivially
parallel $\bO(n)$-cost approximation to the full negative log-likelihood
$\ell(\by) = ( \log \abs{\bS} + \by^T \bS^{-1} \by)/2$, whose naive evaluation
would require~$\bO(n^3)$ for dense $\bS \in \R^{n \times n}$. Even more
conveniently, the approximation (\ref{eq:vecchia_sumform}) directly corresponds
to a matrix-space approximation to the inverse Cholesky factor of $\bS$. Letting
$\bS^{-1} = \bm{T} \bm{T}^T$ be the exact inverse Cholesky for an upper or lower
triangular matrix $\bm{T}$, \cite{pourahmadi1999} describes an interpretation
of~$\bm{T}$ for which its rows or columns correspond to scaled Kriging weights
like $\bm{\lambda}_j$ above, with scaling based on the inverse square root of
the conditional variance $V[y_j \sv y_1, ..., y_{j-1}]$. Considering that the
Vecchia approximation sparsifies these conditional distributions by conditioning
on only a small set of prior enumerated points instead of all of them, this
naturally introduces a sparse approximation $\tilde{\bm{T}} \approx \bm{T}$ with
only  $\bO(n)$ total nonzeros. This connection and construction is summarized in
Figure \ref{fig:vecchia_chol} in the case where $\bm{T}$ is chosen to be upper
triangular (and denoted as $\bm{U}$), and for more details on construction and
derivation we refer the reader to \cite{kaporin1994,katzfuss2021}.
\begin{figure}[!t]
\centering
\scalebox{0.8}{
\begin{tikzpicture}
    \node (left_top) at (0, 1.5) {
        $\begin{aligned}
            &Y_j \mid Y_1, \dots, Y_{j-1} \\
            &\sim \mathcal{N}\left(\sum_{l=1}^{j-1} \lambda_l Y_l, c^2\right)
        \end{aligned}$
    };
    \node (left_top) at (0, 0) {
    $\approx_D$
    };
    \node (left_bot) at (0, -1.5) {
        $\begin{aligned}
            &Y_j \mid Y_{\sigma(j)_1}, \dots, Y_{\sigma(j)_k} \\
            &\sim \mathcal{N}\left(\sum_{l \in \sigma(j)} \tilde{\lambda}_l Y_l, \tilde{c}^2\right)
        \end{aligned}$
    };
    \draw[->, semithick] (2.5, 0) -- (6.5, 0) 
        node[midway, above, yshift=2mm] {$\Sigma^{-1} = UU^T \approx \tilde{U}\tilde{U}^T$};
    \node (right_matrices) at (10.5, 0) {
        $U_{:j} = 
        \left[
        \begin{array}{c}
            -\frac{1}{c}\lambda_1 \\
            -\frac{1}{c}\lambda_2 \\
            -\frac{1}{c}\lambda_3 \\
            \vdots \\
            -\frac{1}{c}\lambda_{j-3} \\
            -\frac{1}{c}\lambda_{j-2} \\
            -\frac{1}{c}\lambda_{j-1} \\
            \hline
            \vspace{-3mm} \\ 
            \frac{1}{c} \\
            \vspace{-3mm} \\
            \hline
            0 \\
            \vdots \\
            0
        \end{array}
        \right],
        \;
        \;
        \tilde{U}_{:j} = 
        \left[
        \begin{array}{c}
            0 \\
            \vdots \\
            -\frac{1}{\tilde{c}}\tilde{\lambda}_{\sigma(j)_1} \\
            0 \\
            \vdots \\
            -\frac{1}{\tilde{c}}\tilde{\lambda}_{\sigma(j)_2} \\
            \vdots \\
            \hline
            \vspace{-3mm} \\
            \frac{1}{\tilde{c}} \\
            \vspace{-3mm} \\
            \hline
            0 \\
            \vdots \\
            0
        \end{array}
        \right]$
    };
\end{tikzpicture}
}
\caption{
A graphical illustration of how the small-sum form of Vecchia approximations
(\ref{eq:vecchia_sumform}) corresponds to a sparse approximation to the inverse
Cholesky factor of $\bS$ with precisely $|\sigma(j)|$ nonzero off-diagonal
entries in each column of $\tilde{\bm{U}}$, using the form $\bS^{-1} = \bm{U}
\bm{U}^T$ (the ``reverse" Cholesky factor) given in \cite{katzfuss2021}.
} 
\label{fig:vecchia_chol}
\end{figure}

With this tool, we can assemble a sparse
preconditioner~$\tilde{\bm{T}} \tilde{\bm{T}}^T$ in $\bO(n)$ time and with
embarrassingly parallel design. Since we will be solving the augmented
  system to obtain weights and polynomial coefficients using PCG, the relevant
  property to study here is the condition number of the preconditioned
  system.
  Define~$\bm{P}_{\text{exact}} = \bm{Q}_{\perp}^T \bS^{-1} \bm{Q}_{\perp}$ to
  be the preconditioner using the exact Mat\'ern matrix $\bS$, and
  $\bm{P}_{\text{Vecchia}} = \bm{Q}_{\perp}^T \tilde{\bm{T}} \tilde{\bm{T}}^T
  \bm{Q}_{\perp}$ to be the preconditioner built with the Vecchia-approximated
  kernel. Letting~$\bm{M}_{\perp} = \bm{Q}_{\perp}^T \bm{M} \bm{Q}_{\perp}$,
  standard matrix norm properties give that
  \begin{equation} 
    \begin{aligned}
  \kappa(\bm{M}_{\perp} \bm{P}_{\text{Vecchia}})
  &\leq
  \kappa(\bm{M}_{\perp} \bm{P}_{\text{exact}}) \, 
  \kappa(\bm{P}_{\text{exact}}^{-1} \bm{P}_{\text{Vecchia}})\\
  &\leq
  \kappa(\bm{M}_{\perp} \bm{P}_{\text{exact}})
    \, \kappa(\bm{T}^{-1} \tilde{\bm{T}})^2.
    \end{aligned}
\end{equation}
Several comments on this bound are in order. First, we note that the bound
\begin{equation}
\kappa(\bm{P}_{\text{exact}}^{-1} \bm{P}_{\text{Vecchia}}) \leq
\kappa(\bm{T}^{-1} \tilde{\bm{T}})^2
\end{equation}
is not a big sacrifice since the Vecchia approximation~$\tilde{\bm{T}}
\tilde{\bm{T}}^{T}$ is constructed so as  to approximate~$\bS^{-1}$ on all
scales, and therefore extending from $(n - S) \times (n - S)$ matrices to the
full $n \times n$ ones is unlikely to have a material impact on a condition
number that should likely be very close to one. However,  the same logic does
not apply to~$\kappa(\bm{M}_{\perp} \bm{P}_{\text{exact}})$ as the inner
matrices~$\bm{M}$ and~$\bS$ \emph{will} materially disagree in their leading
polynomial terms. Conjugating by~$\bm{Q}_{\perp}$ is necessary in order to make
the leading terms in the kernel match. Consulting~\eqref{eq:matern_pre}, we see
that selecting a smaller~$\alpha$ will make the approximations
of~$\bm{M}_{\perp}$ and~$\bm{P}_{\text{exact}}$ sharper, reducing the first
condition number. With that said, however, there is a tradeoff in that
shrinking~$\alpha$ to a sufficiently small value may have a significant impact
on~$\kappa(\bm{T}^{-1} \tilde{\bm{T}})^2$, as the Vecchia approximation may
become less accurate as the process becomes more strongly dependent.

With these observations in mind, we note then that there are two connected
quantities that can be tuned to optimize the efficiency of the PCG solver. The
first is size of the domain, as~\eqref{eq:M_rho} demonstrates that rescaling the
locations of the points by a factor~$\rho$ will have no impact on the final
interpolant. The second quantity is the selection of~$\alpha$ in the Mat\'ern
preconditioner. Intuitively, the balance is between making the leading error
term in (\ref{eq:matern_pre}) smaller, which means a smaller $\rho$ and
$\alpha$, and in making the Vecchia approximation more accurate, which is most
easily achieved by making nearby points less dependent. At present, we do not
offer a theoretically optimal choice of $\rho$ and~$\alpha$, and we note that
the same heuristic motivation given above will apply to any PHS kernel order
$k$, but may have very different optimal choices of $\rho$ and $\alpha$. Using a
gradient-free numerical optimizer in the TPS case $k=2$, scaling locations to be
contained in $[-1/4, 1/4]^2$ and selecting $\alpha \approx 0.35$ produces very
favorable results, often leading to convergence in fewer than $10$ iterations
(as is demonstrated in the next section). For $k=4$, the same procedure suggests
scaling points to be contained in $[-1, 1]^2$ and selecting $\alpha \approx 2$. We
have not investigated higher spline orders, but it is plausible that a careful
analysis would reveal good choices of domain scaling and $\alpha$ for any $k$.

\section{Numerical Examples} \label{sec:demo}

This section contains various numerical experiments to demonstrate the
performance of the scheme described in this paper. Throughout this section, all
Vecchia preconditioners will be assembled with $|\sigma(j)| = 30$ conditioning
points (and full conditioning on prior enumerated points for $j \leq 30$) in the
case of thin-plate splines and $|\sigma(j)| = 50$ for higher-order PHS kernels,
and use a random permutation for the point enumeration. The software companion
to this work, available at
\texttt{https://github.com/cgeoga/FastPolyharmonicSplines.jl}, contains all code
to reproduce the results below.

\subsection{Preconditioning}

To demonstrate the effectiveness of the preconditioner, we consider solving the
system
\begin{equation} \label{eq:pre_compare}
  \bm{M} \bm{c} = \bm{b},
\end{equation}
where $\bm{M}_{j,k} = \tps(\norm[2]{\bx_j - \bx_k})$ is the thin-plate spline
kernel matrix that applied in linear time using the FMM and approach
described in Section~\ref{sec:method}. The locations~$\bx_j$ are taken to be i.i.d.
uniform on $[0,1]^2$, and the right hand side vector~$\bm{b}$ is random with
$\bm{b}_j \sim \Nd(0,1)$.  In previous work~\cite{beatson1999}, an \emph{approximate
cardinal function} preconditioner was proposed which works in a way that is in
some sense very reminiscent of a Vecchia approximation: at each location
$\bx_j$, a small fixed number of nearest neighbors are found (plus a few points
at the endpoints, see~\cite{beatson1999} for details) and local interpolation weights
are computed to build a function that evaluates to one at $\bx_j$ and zero at
all of the other points in the neighborhood set. These weights are assembled
into columns of a sparse matrix, resembling a sparse approximation to the
inverse of $\bm{M}$. Unlike the Vecchia approach we propose here, this sparse
approximate inverse method works directly in the PHS space and thus does not
yield a positive-definite preconditioner. And since it is a direct
preconditioner for $\bm{M}$, for fairness we solve (\ref{eq:pre_compare}) with
GMRES so that the benchmarked approaches are identical save for the
preconditioner design.
\begin{figure}[!t]
  \centering
\begingroup
  \makeatletter
  \providecommand\color[2][]{%
    \GenericError{(gnuplot) \space\space\space\@spaces}{%
      Package color not loaded in conjunction with
      terminal option `colourtext'%
    }{See the gnuplot documentation for explanation.%
    }{Either use 'blacktext' in gnuplot or load the package
      color.sty in LaTeX.}%
    \renewcommand\color[2][]{}%
  }%
  \providecommand\includegraphics[2][]{%
    \GenericError{(gnuplot) \space\space\space\@spaces}{%
      Package graphicx or graphics not loaded%
    }{See the gnuplot documentation for explanation.%
    }{The gnuplot epslatex terminal needs graphicx.sty or graphics.sty.}%
    \renewcommand\includegraphics[2][]{}%
  }%
  \providecommand\rotatebox[2]{#2}%
  \@ifundefined{ifGPcolor}{%
    \newif\ifGPcolor
    \GPcolortrue
  }{}%
  \@ifundefined{ifGPblacktext}{%
    \newif\ifGPblacktext
    \GPblacktexttrue
  }{}%
  \let\gplgaddtomacro\g@addto@macro
  \gdef\gplbacktext{}%
  \gdef\gplfronttext{}%
  \makeatother
  \ifGPblacktext
    \def\colorrgb#1{}%
    \def\colorgray#1{}%
  \else
    \ifGPcolor
      \def\colorrgb#1{\color[rgb]{#1}}%
      \def\colorgray#1{\color[gray]{#1}}%
      \expandafter\def\csname LTw\endcsname{\color{white}}%
      \expandafter\def\csname LTb\endcsname{\color{black}}%
      \expandafter\def\csname LTa\endcsname{\color{black}}%
      \expandafter\def\csname LT0\endcsname{\color[rgb]{1,0,0}}%
      \expandafter\def\csname LT1\endcsname{\color[rgb]{0,1,0}}%
      \expandafter\def\csname LT2\endcsname{\color[rgb]{0,0,1}}%
      \expandafter\def\csname LT3\endcsname{\color[rgb]{1,0,1}}%
      \expandafter\def\csname LT4\endcsname{\color[rgb]{0,1,1}}%
      \expandafter\def\csname LT5\endcsname{\color[rgb]{1,1,0}}%
      \expandafter\def\csname LT6\endcsname{\color[rgb]{0,0,0}}%
      \expandafter\def\csname LT7\endcsname{\color[rgb]{1,0.3,0}}%
      \expandafter\def\csname LT8\endcsname{\color[rgb]{0.5,0.5,0.5}}%
    \else
      \def\colorrgb#1{\color{black}}%
      \def\colorgray#1{\color[gray]{#1}}%
      \expandafter\def\csname LTw\endcsname{\color{white}}%
      \expandafter\def\csname LTb\endcsname{\color{black}}%
      \expandafter\def\csname LTa\endcsname{\color{black}}%
      \expandafter\def\csname LT0\endcsname{\color{black}}%
      \expandafter\def\csname LT1\endcsname{\color{black}}%
      \expandafter\def\csname LT2\endcsname{\color{black}}%
      \expandafter\def\csname LT3\endcsname{\color{black}}%
      \expandafter\def\csname LT4\endcsname{\color{black}}%
      \expandafter\def\csname LT5\endcsname{\color{black}}%
      \expandafter\def\csname LT6\endcsname{\color{black}}%
      \expandafter\def\csname LT7\endcsname{\color{black}}%
      \expandafter\def\csname LT8\endcsname{\color{black}}%
    \fi
  \fi
    \setlength{\unitlength}{0.0500bp}%
    \ifx\gptboxheight\undefined%
      \newlength{\gptboxheight}%
      \newlength{\gptboxwidth}%
      \newsavebox{\gptboxtext}%
    \fi%
    \setlength{\fboxrule}{0.5pt}%
    \setlength{\fboxsep}{1pt}%
    \definecolor{tbcol}{rgb}{1,1,1}%
\begin{picture}(9060.00,3400.00)%
    \gplgaddtomacro\gplbacktext{%
      \csname LTb\endcsname
      \put(803,1189){\makebox(0,0)[r]{\strut{}\footnotesize $10^{-1}$}}%
      \csname LTb\endcsname
      \put(803,1718){\makebox(0,0)[r]{\strut{}\footnotesize $10^{0}$}}%
      \csname LTb\endcsname
      \put(803,2247){\makebox(0,0)[r]{\strut{}\footnotesize $10^{1}$}}%
      \csname LTb\endcsname
      \put(803,2776){\makebox(0,0)[r]{\strut{}\footnotesize $10^{2}$}}%
      \csname LTb\endcsname
      \put(1058,520){\makebox(0,0){\strut{}\footnotesize $2^{10}$}}%
      \csname LTb\endcsname
      \put(1492,520){\makebox(0,0){\strut{}\footnotesize $2^{11}$}}%
      \csname LTb\endcsname
      \put(1926,520){\makebox(0,0){\strut{}\footnotesize $2^{12}$}}%
      \csname LTb\endcsname
      \put(2360,520){\makebox(0,0){\strut{}\footnotesize $2^{13}$}}%
      \csname LTb\endcsname
      \put(2794,520){\makebox(0,0){\strut{}\footnotesize $2^{14}$}}%
      \csname LTb\endcsname
      \put(3228,520){\makebox(0,0){\strut{}\footnotesize $2^{15}$}}%
      \csname LTb\endcsname
      \put(3662,520){\makebox(0,0){\strut{}\footnotesize $2^{16}$}}%
      \csname LTb\endcsname
      \put(4096,520){\makebox(0,0){\strut{}\footnotesize $2^{17}$}}%
    }%
    \gplgaddtomacro\gplfronttext{%
      \csname LTb\endcsname
      \put(1810,2826){\makebox(0,0)[r]{\strut{}\footnotesize $\mathcal{O}(n)$}}%
      \csname LTb\endcsname
      \put(1810,2586){\makebox(0,0)[r]{\strut{}\footnotesize $\mathcal{O}(n^2)$}}%
      \csname LTb\endcsname
      \put(2542,161){\makebox(0,0){\strut{}$n$}}%
      \csname LTb\endcsname
      \put(2542,3221){\makebox(0,0){\strut{}runtime (s)}}%
    }%
    \gplgaddtomacro\gplbacktext{%
      \csname LTb\endcsname
      \put(4758,1059){\makebox(0,0)[r]{\strut{}\footnotesize $10^{1}$}}%
      \csname LTb\endcsname
      \put(4758,2050){\makebox(0,0)[r]{\strut{}\footnotesize $10^{2}$}}%
      \csname LTb\endcsname
      \put(4758,3041){\makebox(0,0)[r]{\strut{}\footnotesize $10^{3}$}}%
      \csname LTb\endcsname
      \put(5013,520){\makebox(0,0){\strut{}\footnotesize $2^{10}$}}%
      \csname LTb\endcsname
      \put(5447,520){\makebox(0,0){\strut{}\footnotesize $2^{11}$}}%
      \csname LTb\endcsname
      \put(5881,520){\makebox(0,0){\strut{}\footnotesize $2^{12}$}}%
      \csname LTb\endcsname
      \put(6315,520){\makebox(0,0){\strut{}\footnotesize $2^{13}$}}%
      \csname LTb\endcsname
      \put(6749,520){\makebox(0,0){\strut{}\footnotesize $2^{14}$}}%
      \csname LTb\endcsname
      \put(7183,520){\makebox(0,0){\strut{}\footnotesize $2^{15}$}}%
      \csname LTb\endcsname
      \put(7617,520){\makebox(0,0){\strut{}\footnotesize $2^{16}$}}%
      \csname LTb\endcsname
      \put(8051,520){\makebox(0,0){\strut{}\footnotesize $2^{17}$}}%
    }%
    \gplgaddtomacro\gplfronttext{%
      \csname LTb\endcsname
      \put(6772,2826){\makebox(0,0)[r]{\strut{}\footnotesize Vecchia}}%
      \csname LTb\endcsname
      \put(6772,2586){\makebox(0,0)[r]{\strut{}\footnotesize Approx. cardinal}}%
      \csname LTb\endcsname
      \put(6497,161){\makebox(0,0){\strut{}$n$}}%
      \csname LTb\endcsname
      \put(6497,3221){\makebox(0,0){\strut{}iteration count}}%
    }%
    \gplbacktext
    \put(0,0){\includegraphics[width={453.00bp},height={170.00bp}]{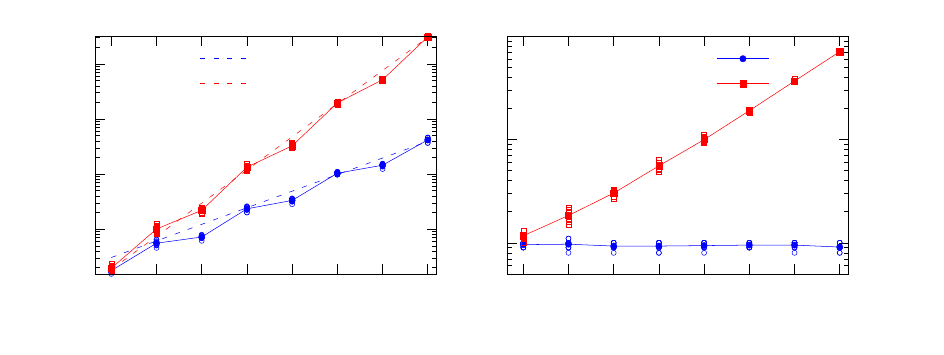}}%
    \gplfronttext
  \end{picture}%
\endgroup

  \caption{A runtime and iteration count comparison for solving the linear
  system (\ref{eq:pre_compare}) using our proposed Vecchia preconditioner with
  $k=50$ conditioning points (blue) and the approximate cardinal function
  preconditioner of \cite{beatson1999} with $50$ neighbors per point (red).
  Solid shapes show the mean of ten trials, and the ten individual runtimes and
  iteration counts are shown as overlaid hollow shapes.}
  \label{fig:runtime}
\end{figure}
Figure \ref{fig:runtime} shows a comparison of the end-to-end runtime and
iteration count of using GMRES to solve (\ref{eq:pre_compare}). We note that in
\cite{beatson1999}, a custom stopping criterion is proposed in an attempt to
better control iteration count. But using standard stopping criteria (in this
case relative and absolute tolerances of $10^{-8}$), we see that iteration
counts are not controlled by the approximate cardinal function preconditioner.
The Vecchia-based preconditioner, on the other hand, averaged fewer than $10$
iterations for all data sizes $n$ despite the fact that the polynomial column
space in which the two kernels differ most was not  removed. We attribute this
improvement to the fact that Vecchia approximations more naturally cover local
and non-local features by virtue of working in square-root space and with a
fixed ordering. Based on the reported iteration counts of \cite{powell1994}, we
expect similar superlinear runtime cost for the preconditioning strategy
reported there as well.

\subsection{Accuracy}

\begin{figure}[!t]
  \centering
\begingroup
  \makeatletter
  \providecommand\color[2][]{%
    \GenericError{(gnuplot) \space\space\space\@spaces}{%
      Package color not loaded in conjunction with
      terminal option `colourtext'%
    }{See the gnuplot documentation for explanation.%
    }{Either use 'blacktext' in gnuplot or load the package
      color.sty in LaTeX.}%
    \renewcommand\color[2][]{}%
  }%
  \providecommand\includegraphics[2][]{%
    \GenericError{(gnuplot) \space\space\space\@spaces}{%
      Package graphicx or graphics not loaded%
    }{See the gnuplot documentation for explanation.%
    }{The gnuplot epslatex terminal needs graphicx.sty or graphics.sty.}%
    \renewcommand\includegraphics[2][]{}%
  }%
  \providecommand\rotatebox[2]{#2}%
  \@ifundefined{ifGPcolor}{%
    \newif\ifGPcolor
    \GPcolortrue
  }{}%
  \@ifundefined{ifGPblacktext}{%
    \newif\ifGPblacktext
    \GPblacktexttrue
  }{}%
  \let\gplgaddtomacro\g@addto@macro
  \gdef\gplbacktext{}%
  \gdef\gplfronttext{}%
  \makeatother
  \ifGPblacktext
    \def\colorrgb#1{}%
    \def\colorgray#1{}%
  \else
    \ifGPcolor
      \def\colorrgb#1{\color[rgb]{#1}}%
      \def\colorgray#1{\color[gray]{#1}}%
      \expandafter\def\csname LTw\endcsname{\color{white}}%
      \expandafter\def\csname LTb\endcsname{\color{black}}%
      \expandafter\def\csname LTa\endcsname{\color{black}}%
      \expandafter\def\csname LT0\endcsname{\color[rgb]{1,0,0}}%
      \expandafter\def\csname LT1\endcsname{\color[rgb]{0,1,0}}%
      \expandafter\def\csname LT2\endcsname{\color[rgb]{0,0,1}}%
      \expandafter\def\csname LT3\endcsname{\color[rgb]{1,0,1}}%
      \expandafter\def\csname LT4\endcsname{\color[rgb]{0,1,1}}%
      \expandafter\def\csname LT5\endcsname{\color[rgb]{1,1,0}}%
      \expandafter\def\csname LT6\endcsname{\color[rgb]{0,0,0}}%
      \expandafter\def\csname LT7\endcsname{\color[rgb]{1,0.3,0}}%
      \expandafter\def\csname LT8\endcsname{\color[rgb]{0.5,0.5,0.5}}%
    \else
      \def\colorrgb#1{\color{black}}%
      \def\colorgray#1{\color[gray]{#1}}%
      \expandafter\def\csname LTw\endcsname{\color{white}}%
      \expandafter\def\csname LTb\endcsname{\color{black}}%
      \expandafter\def\csname LTa\endcsname{\color{black}}%
      \expandafter\def\csname LT0\endcsname{\color{black}}%
      \expandafter\def\csname LT1\endcsname{\color{black}}%
      \expandafter\def\csname LT2\endcsname{\color{black}}%
      \expandafter\def\csname LT3\endcsname{\color{black}}%
      \expandafter\def\csname LT4\endcsname{\color{black}}%
      \expandafter\def\csname LT5\endcsname{\color{black}}%
      \expandafter\def\csname LT6\endcsname{\color{black}}%
      \expandafter\def\csname LT7\endcsname{\color{black}}%
      \expandafter\def\csname LT8\endcsname{\color{black}}%
    \fi
  \fi
    \setlength{\unitlength}{0.0500bp}%
    \ifx\gptboxheight\undefined%
      \newlength{\gptboxheight}%
      \newlength{\gptboxwidth}%
      \newsavebox{\gptboxtext}%
    \fi%
    \setlength{\fboxrule}{0.5pt}%
    \setlength{\fboxsep}{1pt}%
    \definecolor{tbcol}{rgb}{1,1,1}%
\begin{picture}(8500.00,2820.00)%
    \gplgaddtomacro\gplbacktext{%
      \csname LTb\endcsname
      \put(747,792){\makebox(0,0)[r]{\strut{}\small $10^{-1}$}}%
      \csname LTb\endcsname
      \put(747,1232){\makebox(0,0)[r]{\strut{}\small $10^{0}$}}%
      \csname LTb\endcsname
      \put(747,1671){\makebox(0,0)[r]{\strut{}\small $10^{1}$}}%
      \csname LTb\endcsname
      \put(747,2111){\makebox(0,0)[r]{\strut{}\small $10^{2}$}}%
      \csname LTb\endcsname
      \put(848,320){\makebox(0,0){\strut{}\small $2^{8}$}}%
      \csname LTb\endcsname
      \put(1378,320){\makebox(0,0){\strut{}\small $2^{10}$}}%
      \csname LTb\endcsname
      \put(1908,320){\makebox(0,0){\strut{}\small $2^{12}$}}%
      \csname LTb\endcsname
      \put(2438,320){\makebox(0,0){\strut{}\small $2^{14}$}}%
      \csname LTb\endcsname
      \put(2967,320){\makebox(0,0){\strut{}\small $2^{16}$}}%
      \csname LTb\endcsname
      \put(3497,320){\makebox(0,0){\strut{}\small $2^{18}$}}%
      \csname LTb\endcsname
      \put(4027,320){\makebox(0,0){\strut{}\small $2^{20}$}}%
    }%
    \gplgaddtomacro\gplfronttext{%
      \csname LTb\endcsname
      \put(2437,44){\makebox(0,0){\strut{}\small $n$}}%
      \csname LTb\endcsname
      \put(2437,2687){\makebox(0,0){\strut{}\small runtime (s)}}%
    }%
    \gplgaddtomacro\gplbacktext{%
      \csname LTb\endcsname
      \put(4775,762){\makebox(0,0)[r]{\strut{}\small $10^{-6}$}}%
      \csname LTb\endcsname
      \put(4775,1052){\makebox(0,0)[r]{\strut{}\small $10^{-5}$}}%
      \csname LTb\endcsname
      \put(4775,1341){\makebox(0,0)[r]{\strut{}\small $10^{-4}$}}%
      \csname LTb\endcsname
      \put(4775,1631){\makebox(0,0)[r]{\strut{}\small $10^{-3}$}}%
      \csname LTb\endcsname
      \put(4775,1921){\makebox(0,0)[r]{\strut{}\small $10^{-2}$}}%
      \csname LTb\endcsname
      \put(4775,2210){\makebox(0,0)[r]{\strut{}\small $10^{-1}$}}%
      \csname LTb\endcsname
      \put(4775,2500){\makebox(0,0)[r]{\strut{}\small $10^{0}$}}%
      \csname LTb\endcsname
      \put(4875,320){\makebox(0,0){\strut{}\small $2^{8}$}}%
      \csname LTb\endcsname
      \put(5405,320){\makebox(0,0){\strut{}\small $2^{10}$}}%
      \csname LTb\endcsname
      \put(5935,320){\makebox(0,0){\strut{}\small $2^{12}$}}%
      \csname LTb\endcsname
      \put(6465,320){\makebox(0,0){\strut{}\small $2^{14}$}}%
      \csname LTb\endcsname
      \put(6995,320){\makebox(0,0){\strut{}\small $2^{16}$}}%
      \csname LTb\endcsname
      \put(7525,320){\makebox(0,0){\strut{}\small $2^{18}$}}%
      \csname LTb\endcsname
      \put(8055,320){\makebox(0,0){\strut{}\small $2^{20}$}}%
    }%
    \gplgaddtomacro\gplfronttext{%
      \csname LTb\endcsname
      \put(6560,1322){\makebox(0,0)[r]{\strut{}\footnotesize exact}}%
      \csname LTb\endcsname
      \put(6560,1082){\makebox(0,0)[r]{\strut{}\footnotesize approx. (tol. $10^{-5}$)}}%
      \csname LTb\endcsname
      \put(6560,842){\makebox(0,0)[r]{\strut{}\footnotesize approx. (tol. $10^{-10}$)}}%
      \csname LTb\endcsname
      \put(6465,44){\makebox(0,0){\strut{}\small $n$}}%
      \csname LTb\endcsname
      \put(6465,2687){\makebox(0,0){\strut{}\small $\ell_{\infty}$ prediction error}}%
    }%
    \gplbacktext
    \put(0,0){\includegraphics[width={425.00bp},height={141.00bp}]{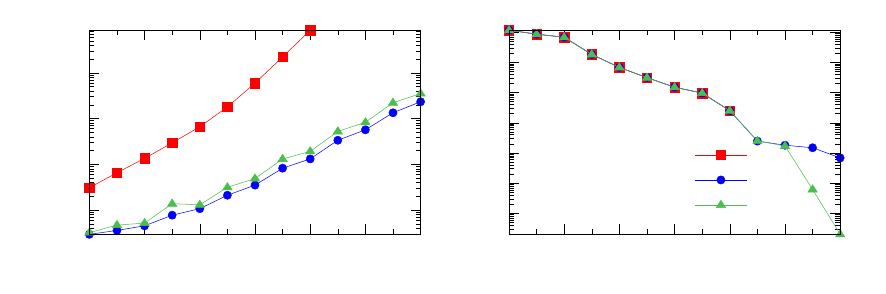}}%
    \gplfronttext
  \end{picture}%
\endgroup

  \caption{
  Runtime and prediction accuracy verifications for the problem of predicting
  $f(x)$ given by (\ref{eq:f_interp}) on a $100 \times 100$ regular grid on
  $[0,1]^2$ using $n$ points from a Sobol sequence. For the blue/circle line,
  the FMM was configured to use a tolerance of $\varepsilon = 10^{-8}$ and PCG
  configured to declare convergence at relative and absolute tolerances of
  $10^{-5}$. For the triangle/green line, both the FMM and PCG tolerances were
  set to $10^{-10}$.
  }
  \label{fig:errortime}
\end{figure}

We now assess the accuracy in computing the weights~$\bm{w}$ using the
accelerated application of $\bv \mapsto
\bm{M} \bv$ and the PCG solve. As a baseline, we sample the function
\begin{equation} \label{eq:f_interp}
  f(x_1, x_2) = \cos(30 x_1 + 0.8) \cdot \sin(30 x_2 - 2.1) 
  + 
  5 \exp\left( -5 \sqrt{(x_1 - 0.25)^2 + (x_2 - 0.75)^2} \right)
\end{equation}
and make predictions at densely gridded locations using~\eqref{eq:f_interp}.
This choice of function poses a challenge for global interpolators by requiring
the weights to resolve both a smooth oscillatory function and a $\mathcal{C}^0$
singularity at one point. Moreover, this singular point ensures that the PCG
solver for weights $\bm{w}$ will need to contend with small eigenvalues of the
matrix $\Qperp^T \bm{M} \Qperp$ in the solution of~\eqref{eq:phs_linsys_solve},
therefore also serving as a stress test for the preconditioner. For these
reasons, it is suitable for a large-$n$ performance analysis. Figure
\ref{fig:errortime} shows results from the simulation study in which $n$ points
are given from a Sobol sequence \cite{sobol1967} on $[0,1]^2$ and then
``regridded", meaning that the interpolating function (\ref{eq:phs_interp}) is
prepared and then evaluated on a dense $100 \times 100$ grid of points on
$[0,1]^2$.

Figure \ref{fig:errortime} reports the end-to-end runtime cost of using dense
linear algebra to exactly solve system~\eqref{eq:phs_linsys} and
evaluate~\eqref{eq:phs_interp} versus the accelerated method proposed here as
well as an error comparison with standard $\bO(n^3)$ dense methods.
Runtimes include the cost of preconditioner assembly. All code was run
using~$16$ cores of an AMD EPYC~$9554$P. Importantly, we see that the
iterative FMM-accelerated PCG routine produces identical errors to the exact
method that uses a Bunch-Kaufman factorization on the augmented matrix
in~\eqref{eq:phs_linsys}, and that at higher tolerances the solvers are
sufficiently exact to show $\ell_{\infty}$ error improvements with over $n=1M$
points. Not shown here is the $\ell_2$ prediction error, as in all cases the
errors agreed to sufficient precision that visualization is not productive.

\subsection{Higher-order interpolation}

\begin{figure}[!t]
  \centering
\begingroup
  \makeatletter
  \providecommand\color[2][]{%
    \GenericError{(gnuplot) \space\space\space\@spaces}{%
      Package color not loaded in conjunction with
      terminal option `colourtext'%
    }{See the gnuplot documentation for explanation.%
    }{Either use 'blacktext' in gnuplot or load the package
      color.sty in LaTeX.}%
    \renewcommand\color[2][]{}%
  }%
  \providecommand\includegraphics[2][]{%
    \GenericError{(gnuplot) \space\space\space\@spaces}{%
      Package graphicx or graphics not loaded%
    }{See the gnuplot documentation for explanation.%
    }{The gnuplot epslatex terminal needs graphicx.sty or graphics.sty.}%
    \renewcommand\includegraphics[2][]{}%
  }%
  \providecommand\rotatebox[2]{#2}%
  \@ifundefined{ifGPcolor}{%
    \newif\ifGPcolor
    \GPcolortrue
  }{}%
  \@ifundefined{ifGPblacktext}{%
    \newif\ifGPblacktext
    \GPblacktexttrue
  }{}%
  \let\gplgaddtomacro\g@addto@macro
  \gdef\gplbacktext{}%
  \gdef\gplfronttext{}%
  \makeatother
  \ifGPblacktext
    \def\colorrgb#1{}%
    \def\colorgray#1{}%
  \else
    \ifGPcolor
      \def\colorrgb#1{\color[rgb]{#1}}%
      \def\colorgray#1{\color[gray]{#1}}%
      \expandafter\def\csname LTw\endcsname{\color{white}}%
      \expandafter\def\csname LTb\endcsname{\color{black}}%
      \expandafter\def\csname LTa\endcsname{\color{black}}%
      \expandafter\def\csname LT0\endcsname{\color[rgb]{1,0,0}}%
      \expandafter\def\csname LT1\endcsname{\color[rgb]{0,1,0}}%
      \expandafter\def\csname LT2\endcsname{\color[rgb]{0,0,1}}%
      \expandafter\def\csname LT3\endcsname{\color[rgb]{1,0,1}}%
      \expandafter\def\csname LT4\endcsname{\color[rgb]{0,1,1}}%
      \expandafter\def\csname LT5\endcsname{\color[rgb]{1,1,0}}%
      \expandafter\def\csname LT6\endcsname{\color[rgb]{0,0,0}}%
      \expandafter\def\csname LT7\endcsname{\color[rgb]{1,0.3,0}}%
      \expandafter\def\csname LT8\endcsname{\color[rgb]{0.5,0.5,0.5}}%
    \else
      \def\colorrgb#1{\color{black}}%
      \def\colorgray#1{\color[gray]{#1}}%
      \expandafter\def\csname LTw\endcsname{\color{white}}%
      \expandafter\def\csname LTb\endcsname{\color{black}}%
      \expandafter\def\csname LTa\endcsname{\color{black}}%
      \expandafter\def\csname LT0\endcsname{\color{black}}%
      \expandafter\def\csname LT1\endcsname{\color{black}}%
      \expandafter\def\csname LT2\endcsname{\color{black}}%
      \expandafter\def\csname LT3\endcsname{\color{black}}%
      \expandafter\def\csname LT4\endcsname{\color{black}}%
      \expandafter\def\csname LT5\endcsname{\color{black}}%
      \expandafter\def\csname LT6\endcsname{\color{black}}%
      \expandafter\def\csname LT7\endcsname{\color{black}}%
      \expandafter\def\csname LT8\endcsname{\color{black}}%
    \fi
  \fi
    \setlength{\unitlength}{0.0500bp}%
    \ifx\gptboxheight\undefined%
      \newlength{\gptboxheight}%
      \newlength{\gptboxwidth}%
      \newsavebox{\gptboxtext}%
    \fi%
    \setlength{\fboxrule}{0.5pt}%
    \setlength{\fboxsep}{1pt}%
    \definecolor{tbcol}{rgb}{1,1,1}%
\begin{picture}(8500.00,2820.00)%
    \gplgaddtomacro\gplbacktext{%
      \csname LTb\endcsname
      \put(747,781){\makebox(0,0)[r]{\strut{}\small $10^{-1}$}}%
      \csname LTb\endcsname
      \put(747,1223){\makebox(0,0)[r]{\strut{}\small $10^{0}$}}%
      \csname LTb\endcsname
      \put(747,1665){\makebox(0,0)[r]{\strut{}\small $10^{1}$}}%
      \csname LTb\endcsname
      \put(747,2107){\makebox(0,0)[r]{\strut{}\small $10^{2}$}}%
      \csname LTb\endcsname
      \put(1197,320){\makebox(0,0){\strut{}\small $2^{10}$}}%
      \csname LTb\endcsname
      \put(1550,320){\makebox(0,0){\strut{}\small $2^{12}$}}%
      \csname LTb\endcsname
      \put(1904,320){\makebox(0,0){\strut{}\small $2^{14}$}}%
      \csname LTb\endcsname
      \put(2257,320){\makebox(0,0){\strut{}\small $2^{16}$}}%
    }%
    \gplgaddtomacro\gplfronttext{%
      \csname LTb\endcsname
      \put(1554,44){\makebox(0,0){\strut{}\small $n$}}%
      \csname LTb\endcsname
      \put(1554,2687){\makebox(0,0){\strut{}\small runtime (s)}}%
    }%
    \gplgaddtomacro\gplbacktext{%
      \csname LTb\endcsname
      \put(3008,714){\makebox(0,0)[r]{\strut{}\small $10^{-6}$}}%
      \csname LTb\endcsname
      \put(3008,1097){\makebox(0,0)[r]{\strut{}\small $10^{-5}$}}%
      \csname LTb\endcsname
      \put(3008,1481){\makebox(0,0)[r]{\strut{}\small $10^{-4}$}}%
      \csname LTb\endcsname
      \put(3008,1864){\makebox(0,0)[r]{\strut{}\small $10^{-3}$}}%
      \csname LTb\endcsname
      \put(3008,2247){\makebox(0,0)[r]{\strut{}\small $10^{-2}$}}%
      \csname LTb\endcsname
      \put(3458,320){\makebox(0,0){\strut{}\small $2^{10}$}}%
      \csname LTb\endcsname
      \put(3811,320){\makebox(0,0){\strut{}\small $2^{12}$}}%
      \csname LTb\endcsname
      \put(4165,320){\makebox(0,0){\strut{}\small $2^{14}$}}%
      \csname LTb\endcsname
      \put(4518,320){\makebox(0,0){\strut{}\small $2^{16}$}}%
    }%
    \gplgaddtomacro\gplfronttext{%
      \csname LTb\endcsname
      \put(3815,44){\makebox(0,0){\strut{}\small $n$}}%
      \csname LTb\endcsname
      \put(3815,2687){\makebox(0,0){\strut{}\small $\ell_{\infty}$ error}}%
    }%
    \gplgaddtomacro\gplbacktext{%
      \csname LTb\endcsname
      \put(5269,643){\makebox(0,0)[r]{\strut{}\small $10^{-6}$}}%
      \csname LTb\endcsname
      \put(5269,1002){\makebox(0,0)[r]{\strut{}\small $10^{-5}$}}%
      \csname LTb\endcsname
      \put(5269,1361){\makebox(0,0)[r]{\strut{}\small $10^{-4}$}}%
      \csname LTb\endcsname
      \put(5269,1721){\makebox(0,0)[r]{\strut{}\small $10^{-3}$}}%
      \csname LTb\endcsname
      \put(5269,2080){\makebox(0,0)[r]{\strut{}\small $10^{-2}$}}%
      \csname LTb\endcsname
      \put(5269,2439){\makebox(0,0)[r]{\strut{}\small $10^{-1}$}}%
      \csname LTb\endcsname
      \put(5719,320){\makebox(0,0){\strut{}\small $2^{10}$}}%
      \csname LTb\endcsname
      \put(6073,320){\makebox(0,0){\strut{}\small $2^{12}$}}%
      \csname LTb\endcsname
      \put(6426,320){\makebox(0,0){\strut{}\small $2^{14}$}}%
      \csname LTb\endcsname
      \put(6779,320){\makebox(0,0){\strut{}\small $2^{16}$}}%
    }%
    \gplgaddtomacro\gplfronttext{%
      \csname LTb\endcsname
      \put(7981,2400){\makebox(0,0)[r]{\strut{}\footnotesize exact}}%
      \csname LTb\endcsname
      \put(7981,2160){\makebox(0,0)[r]{\strut{}\footnotesize $k=2$ (FMM)}}%
      \csname LTb\endcsname
      \put(7981,1920){\makebox(0,0)[r]{\strut{}\footnotesize $k=4$ (FMM)}}%
      \csname LTb\endcsname
      \put(7981,1680){\makebox(0,0)[r]{\strut{}\footnotesize $k=4$ (hybrid)}}%
      \csname LTb\endcsname
      \put(6077,44){\makebox(0,0){\strut{}\small $n$}}%
      \csname LTb\endcsname
      \put(6077,2687){\makebox(0,0){\strut{}\small $\ell_{2}$ error}}%
    }%
    \gplbacktext
    \put(0,0){\includegraphics[width={425.00bp},height={141.00bp}]{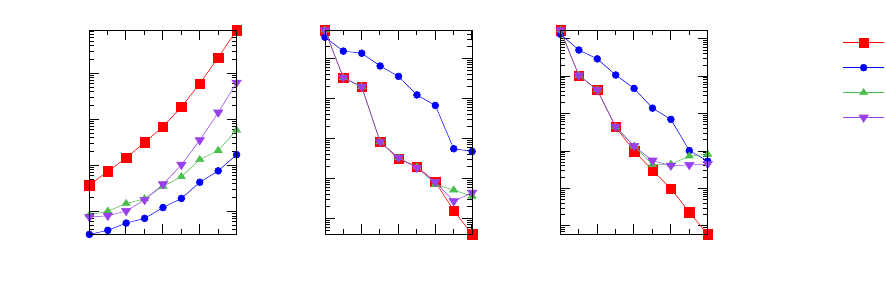}}%
    \gplfronttext
  \end{picture}%
\endgroup

  \caption{
  Runtime and prediction accuracy verifications for the problem of predicting
  the Franke function on a $100 \times 100$ regular grid on $[0,1]^2$ using $n$
  points from a Sobol sequence. As above, the red/square line gives exact
  solutions, and the blue/circle line gives PHS interpolation with $k=2$
  (thin-plate splines). The green/triangle line gives the FMM and
  PCG-accelerated analog with the PHS $\varphi$ with $k=4$, and the purple line
  gives a hybrid routine that uses exact matrix-vector products within the PCG
  framework.
  }
  \label{fig:phsk4_errortime}
\end{figure}

Next, we demonstrate the runtime cost and error of higher-order PHS
interpolation. In particular, we consider the same sampling scenario as above,
but now we produce regridding interpolants for the Franke
function~\cite{franke1982}, which is a small sum of Gaussian and exponential
functions given by
\begin{align*} 
  f_{\text{Franke}}(\bm{x}) = 
  \sum_{j=1}^{4} w_j \exp\left( -(9\bm{x} - \bm{c}_j)^T \bm{\Lambda}_j (9\bm{x}
  - \bm{c}_j) \right),
\end{align*}
where $\bm{x} = [x, y]^T$, weights are $\bm{w} = (0.75, 0.75, 0.50, -0.20)$,
centers are $\bm{c}_1 = (2, 2)$, $\bm{c}_2 = (-1, -1)$, $\bm{c}_3 = (7,
3)$, $\bm{c}_4 = (4, 7)$, and diagonal metric matrices are $\bm{\Lambda}_1 =
\bm{\Lambda}_3 = \frac{1}{4}\mathbf{I}_2$, $\bm{\Lambda}_2 =
\operatorname{diag}\left(\frac{1}{49}, \frac{1}{10}\right)$, and $\bm{\Lambda}_4
= \mathbf{I}_2$.  Since the new function to approximate is analytic everywhere,
it is beneficial to use higher-order PHS functions instead of merely $\tps$.
Figure \ref{fig:phsk4_errortime} demonstrates this phenomenon, and also
validates the preserved runtime cost of the proposed approach for higher-order
PHS functions.  Since even for $k=4$ the condition number of the kernel matrix
$\bm{M}$ can exceed~$10^{16}$ for even moderate $n$, we do not take~$n$ to such
high of values as in Figure~\ref{fig:errortime}. Nonetheless, however, we see
that runtimes still compare favorably with exact methods. For completeness, we
also add a purple/triangle line representing the approach that uses an exact
matrix-vector product for $\bm{M}$, naturally costing $\bO(n^2)$ work and
memory, but still using the PCG solver and Vecchia preconditioner. 

We offer several comments on the obtained results:  First, we see that the clear
linear cost of the end-to-end pipeline is preserved for the PHS with $k=4$,
albeit it with a higher prefactor than with using $\tps$. There are several
reasons for this. Most obviously, applying the PHS matrix $\bm{M}$ with $k=4$
requires $10$ FMMs instead of $4$ for $k=2$, since the rank of the matrix with
entries $\norm[2]{\bx_j - \bx_k}^4$ is higher due to additional cross terms. The
\texttt{fmm2d} library has excellent vectorized routines, and so the cost is
materially less than $10$ individual FMMs, but naturally it is still higher.
Further, for $k=4$, the PHS matrix $\bm{M}$ is significantly more
ill-conditioned than with $k=2$. And while the Vecchia preconditioner still
controls iteration count with $n$, the average number of iterations required was
around $15$ instead of $5$. Finally, we note that due to this poor conditioning
and sensitivity, it is prudent to tighten the tolerance of the FMM algorithm
itself. The results of Figure \ref{fig:phsk4_errortime} were computed with an
FMM tolerance of $\varepsilon = 10^{-12}$ for all results (including $k=2$),
which is why the runtimes for interpolating with $\tps$ are also slightly longer
than they were in Figure \ref{fig:errortime}.

Second, we see that the increasingly poor conditioning of the PHS matrix
$\bm{M}$ inevitably makes it harder for the PCG solver to produce weights that
work exactly as well as the exact dense solver, which uses a numerically stable
Bunch-Kaufman factorization. It is well-known in the literature that
higher-order PHS methods are typically not used globally and are not typically
taken to such significant data sizes. But with that said, as Figure
\ref{fig:phsk4_errortime} demonstrates that for even moderate data sizes the speedup
(and memory savings) that this approach offer are significant.

\subsection{Application to real data}

\begin{figure}[!ht]
  \centering
\begingroup
  \makeatletter
  \providecommand\color[2][]{%
    \GenericError{(gnuplot) \space\space\space\@spaces}{%
      Package color not loaded in conjunction with
      terminal option `colourtext'%
    }{See the gnuplot documentation for explanation.%
    }{Either use 'blacktext' in gnuplot or load the package
      color.sty in LaTeX.}%
    \renewcommand\color[2][]{}%
  }%
  \providecommand\includegraphics[2][]{%
    \GenericError{(gnuplot) \space\space\space\@spaces}{%
      Package graphicx or graphics not loaded%
    }{See the gnuplot documentation for explanation.%
    }{The gnuplot epslatex terminal needs graphicx.sty or graphics.sty.}%
    \renewcommand\includegraphics[2][]{}%
  }%
  \providecommand\rotatebox[2]{#2}%
  \@ifundefined{ifGPcolor}{%
    \newif\ifGPcolor
    \GPcolortrue
  }{}%
  \@ifundefined{ifGPblacktext}{%
    \newif\ifGPblacktext
    \GPblacktexttrue
  }{}%
  \let\gplgaddtomacro\g@addto@macro
  \gdef\gplbacktext{}%
  \gdef\gplfronttext{}%
  \makeatother
  \ifGPblacktext
    \def\colorrgb#1{}%
    \def\colorgray#1{}%
  \else
    \ifGPcolor
      \def\colorrgb#1{\color[rgb]{#1}}%
      \def\colorgray#1{\color[gray]{#1}}%
      \expandafter\def\csname LTw\endcsname{\color{white}}%
      \expandafter\def\csname LTb\endcsname{\color{black}}%
      \expandafter\def\csname LTa\endcsname{\color{black}}%
      \expandafter\def\csname LT0\endcsname{\color[rgb]{1,0,0}}%
      \expandafter\def\csname LT1\endcsname{\color[rgb]{0,1,0}}%
      \expandafter\def\csname LT2\endcsname{\color[rgb]{0,0,1}}%
      \expandafter\def\csname LT3\endcsname{\color[rgb]{1,0,1}}%
      \expandafter\def\csname LT4\endcsname{\color[rgb]{0,1,1}}%
      \expandafter\def\csname LT5\endcsname{\color[rgb]{1,1,0}}%
      \expandafter\def\csname LT6\endcsname{\color[rgb]{0,0,0}}%
      \expandafter\def\csname LT7\endcsname{\color[rgb]{1,0.3,0}}%
      \expandafter\def\csname LT8\endcsname{\color[rgb]{0.5,0.5,0.5}}%
    \else
      \def\colorrgb#1{\color{black}}%
      \def\colorgray#1{\color[gray]{#1}}%
      \expandafter\def\csname LTw\endcsname{\color{white}}%
      \expandafter\def\csname LTb\endcsname{\color{black}}%
      \expandafter\def\csname LTa\endcsname{\color{black}}%
      \expandafter\def\csname LT0\endcsname{\color{black}}%
      \expandafter\def\csname LT1\endcsname{\color{black}}%
      \expandafter\def\csname LT2\endcsname{\color{black}}%
      \expandafter\def\csname LT3\endcsname{\color{black}}%
      \expandafter\def\csname LT4\endcsname{\color{black}}%
      \expandafter\def\csname LT5\endcsname{\color{black}}%
      \expandafter\def\csname LT6\endcsname{\color{black}}%
      \expandafter\def\csname LT7\endcsname{\color{black}}%
      \expandafter\def\csname LT8\endcsname{\color{black}}%
    \fi
  \fi
    \setlength{\unitlength}{0.0500bp}%
    \ifx\gptboxheight\undefined%
      \newlength{\gptboxheight}%
      \newlength{\gptboxwidth}%
      \newsavebox{\gptboxtext}%
    \fi%
    \setlength{\fboxrule}{0.5pt}%
    \setlength{\fboxsep}{1pt}%
    \definecolor{tbcol}{rgb}{1,1,1}%
\begin{picture}(9060.00,2820.00)%
    \gplgaddtomacro\gplbacktext{%
      \csname LTb\endcsname
      \put(803,690){\makebox(0,0)[r]{\strut{}\scriptsize 12}}%
      \csname LTb\endcsname
      \put(803,971){\makebox(0,0)[r]{\strut{}\scriptsize 14}}%
      \csname LTb\endcsname
      \put(803,1252){\makebox(0,0)[r]{\strut{}\scriptsize 16}}%
      \csname LTb\endcsname
      \put(803,1533){\makebox(0,0)[r]{\strut{}\scriptsize 18}}%
      \csname LTb\endcsname
      \put(803,1814){\makebox(0,0)[r]{\strut{}\scriptsize 20}}%
      \csname LTb\endcsname
      \put(803,2095){\makebox(0,0)[r]{\strut{}\scriptsize 22}}%
      \csname LTb\endcsname
      \put(803,2375){\makebox(0,0)[r]{\strut{}\scriptsize 24}}%
      \csname LTb\endcsname
      \put(1110,360){\makebox(0,0){\strut{}\scriptsize -114}}%
      \csname LTb\endcsname
      \put(1914,360){\makebox(0,0){\strut{}\scriptsize -110}}%
      \csname LTb\endcsname
      \put(2718,360){\makebox(0,0){\strut{}\scriptsize -106}}%
      \csname LTb\endcsname
      \put(3522,360){\makebox(0,0){\strut{}\scriptsize -102}}%
    }%
    \gplgaddtomacro\gplfronttext{%
      \csname LTb\endcsname
      \put(406,1469){\rotatebox{-270.00}{\makebox(0,0){\strut{}\scriptsize Latitude}}}%
      \csname LTb\endcsname
      \put(2406,60){\makebox(0,0){\strut{}\scriptsize Longitude}}%
      \csname LTb\endcsname
      \put(2406,2759){\makebox(0,0){\strut{}\scriptsize full SST anomaly data (deg. C)}}%
    }%
    \gplgaddtomacro\gplbacktext{%
      \csname LTb\endcsname
      \put(4577,690){\makebox(0,0)[r]{\strut{}\scriptsize 12}}%
      \csname LTb\endcsname
      \put(4577,971){\makebox(0,0)[r]{\strut{}\scriptsize 14}}%
      \csname LTb\endcsname
      \put(4577,1252){\makebox(0,0)[r]{\strut{}\scriptsize 16}}%
      \csname LTb\endcsname
      \put(4577,1533){\makebox(0,0)[r]{\strut{}\scriptsize 18}}%
      \csname LTb\endcsname
      \put(4577,1814){\makebox(0,0)[r]{\strut{}\scriptsize 20}}%
      \csname LTb\endcsname
      \put(4577,2095){\makebox(0,0)[r]{\strut{}\scriptsize 22}}%
      \csname LTb\endcsname
      \put(4577,2375){\makebox(0,0)[r]{\strut{}\scriptsize 24}}%
      \csname LTb\endcsname
      \put(4884,360){\makebox(0,0){\strut{}\scriptsize -114}}%
      \csname LTb\endcsname
      \put(5688,360){\makebox(0,0){\strut{}\scriptsize -110}}%
      \csname LTb\endcsname
      \put(6492,360){\makebox(0,0){\strut{}\scriptsize -106}}%
      \csname LTb\endcsname
      \put(7296,360){\makebox(0,0){\strut{}\scriptsize -102}}%
    }%
    \gplgaddtomacro\gplfronttext{%
      \csname LTb\endcsname
      \put(6181,60){\makebox(0,0){\strut{}\scriptsize Longitude}}%
      \csname LTb\endcsname
      \put(8009,677){\makebox(0,0)[l]{\strut{}\scriptsize -2}}%
      \csname LTb\endcsname
      \put(8009,1073){\makebox(0,0)[l]{\strut{}\scriptsize -1}}%
      \csname LTb\endcsname
      \put(8009,1469){\makebox(0,0)[l]{\strut{}\scriptsize 0}}%
      \csname LTb\endcsname
      \put(8009,1866){\makebox(0,0)[l]{\strut{}\scriptsize 1}}%
      \csname LTb\endcsname
      \put(8009,2262){\makebox(0,0)[l]{\strut{}\scriptsize 2}}%
      \csname LTb\endcsname
      \put(6181,2759){\makebox(0,0){\strut{}\scriptsize cloud gap SST anomaly data (deg. C)}}%
    }%
    \gplbacktext
    \put(0,0){\includegraphics[width={453.00bp},height={141.00bp}]{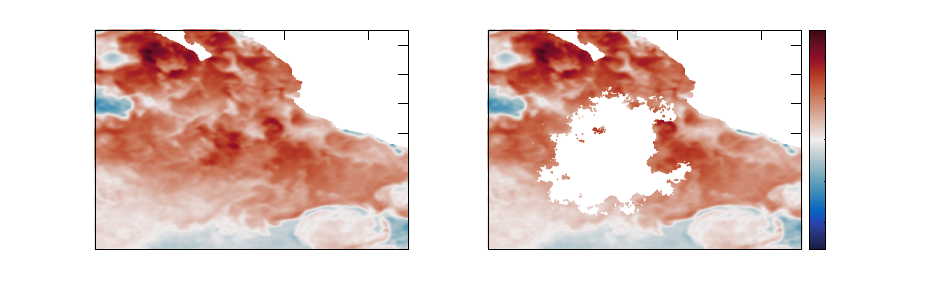}}%
    \gplfronttext
  \end{picture}%
\endgroup

  \caption{The full sea surface temperature anomaly data (left) and synthetic
  cloud-based missingness pattern (right).}
  \label{fig:ssta}
\end{figure}

As a final demonstration, we give an example of interpolating real data. In
particular, we take sea surface temperature anomaly measurements from the NOAA
Coral Reef Watch database \cite{noaa_crw} in the central Pacific ocean and use
the synthetic missingness pattern resembling measurement gaps due to cloud cover
from \cite{arora2026}. This is a common issue in satellite-based remote sensing
datasets, leaving approximately $n=58\,000$ measurements with which to predict
the missing ones.  In a traditional interpolation approach using scalable
Gaussian process methods, one would specify a mean and covariance function (like
the Mat\'ern model that this work uses as a preconditioner), estimate parameters
somehow (for example using a Vecchia approximation), and then perform some
variety of nearest neighbor-based interpolation. Due to the scale-invariance of
polyharmonic splines, however, if one is only interested in prediction and not
also performing uncertainty quantification, PHS interpolants can often perform
very well even with no parameter tuning.

\begin{table}[!ht]
\label{tab:ssta}
\centering
\begin{tabular}{ccc}
Method & runtime (s) & RMSE \\
\hline
Thin-plate spline                       & $\bm{5}$   & $\bm{0.287}$ \\
Mat\'ern + Vecchia ($|\sigma(j)| = 300$)  & $418$      & $0.398$     
\end{tabular}
\caption{A summary of the runtime cost and prediction RMSE for interpolating the
missing cloud-covered measurements shown in the right panel of Figure
\ref{fig:ssta}. Timings were obtained on an ultrabook form factor laptop with in
Intel i7-1260P processor, and in the Mat\'ern model case include parameter
estimation as well as computing predictions.} 
\end{table}

Figure \ref{fig:ssta} shows the full data and the partial data after removing
the measurements under synthetic cloud cover, and Table $1$ shows the runtime
cost and prediction RMSE of thin-plate spline interpolation and Kriging with a
Mat\'ern covariance function and a Vecchia approximation using standard
high-performance design choices. As Table $1$ demonstrates, the thin-plate
spline interpolation performs noticeably better than the Mat\'ern kernel and
Vecchia approximation despite not having any tuning parameters to select. This
disparity in performance highlights an interesting difference between the two
functions: while the Mat\'ern kernel and $\tps$ have similar behavior near the
origin, the pairwise covariances implied by $\bm{Q}_{\perp}^T \bm{M}
\bm{Q}_{\perp}$ and $\bm{Q}_{\perp}^T \bS \bm{Q}_{\perp}$ are actually very
different, exhibiting very slow decay (a model feature known as \emph{long
memory}) in the former case and exponential decay in the latter
\cite{NIST,matheron1973}. Due to the scale invariance of PHS kernels, however,
the Mat\'ern model can nonetheless serve as an extremely effective
preconditioner with suitable point rescalings and parameter matching as
described in Section \ref{subsec:vecchia}.

\section{Discussion} \label{sec:discussion}

This work proposes a method for polyharmonic spline interpolation that combines
several existing high-performance and well-studied methods---the fast multipole
method (FMM) \cite{greengard1987fast} for matrix-vector products and the Vecchia
approximation \cite{vecchia1988,kaporin1994} for preconditioner design---with a
few straightforward strategies for generalizing over the spline order and
enabling the use of PCG for iteratively solving linear systems. Putting these
observations together leads to a general and very performant method for PHS
interpolation with materially better runtime performance than prior approaches
\emph{and} effectively exact agreement with exact dense methods until the linear
system involves a matrix whose spectrum exhausts the dynamic range of double
precision floating point numbers. The software library that is companion to this
work offers a specially tuned and performant preconditioner design for~$k=2$
and~$k=4$.  Producing similarly tuned performant methods for matrix-vector
products with $\bm{M}_{j,k} = \norm[2]{\bx_j - \bx_k}$ will unlock similar
performance for $k \in \set{1,3}$, and is an exciting avenue for future work, as
is extending to three and higher dimensions. 

Several other potentially useful follow-up projects to this work would also be
natural. For one, we note that the kernel $\Gamma(s)^{-1} r^{s}$ for $s > 0$,
(with the same limiting case of above of $\Gamma(s)^{-1} r^{s} \log r$ when $s$
is an even integer), is a more general power law kernel that exhibits similar
scale-invariance properties but provides more sensitivity to local smoothness.
An algorithm for rapid matrix-vector products with that kernel would connect
more directly with the generalized power law models in the Gaussian process
literature \cite{matheron1973}. Similarly, we note that scale-invariance in
predictions is a property that is specific to first-order information. If one
were also interested in obtaining confidence intervals for predictions, for
example, then tuning parameters would be necessary to optimize exactly as one
does for standard GP covariance functions.  For that reason, log-determinants
of~$\bm{M}_{\perp}$ would be necessary. The preconditioning tools here may be
applicable to stochastic log-determinant estimation routines like stochastic
Lanczos quadrature~\cite{ubaru2017}.

\section*{Acknowledgments}
This material is based on work supported by the National Science Foundation
under Grant Number DMS-2610202.  Any opinions, findings, and conclusions or
recommendations expressed in this material are those of the author(s) and do not
necessarily reflect the views of the National Science Foundation.

\bibliographystyle{abbrv}
\bibliography{references}

\end{document}